\documentclass[a4paper,onecolumn,11pt,accepted=0000-00-00]{quantumarticle}
\pdfoutput=1
\usepackage[utf8]{inputenc}
\usepackage[english]{babel}
\usepackage[T1]{fontenc}
\usepackage{amsmath}
\usepackage{hyperref}
\usepackage{orcidlink}

\usepackage{tikz}
\usepackage{lipsum}

\usepackage{booktabs}
\usepackage{multirow}
\usepackage{adjustbox}

\usepackage[backend=bibtex,sorting=none]{biblatex}
\usepackage{float}

\newcommand{\e}{\mathbf{e}}

\newcommand{\st}{\textrm{s.t.}}

\newcommand{\Diag}{\textrm{Diag}}

\newcommand\x{{\mathbf x}}

\begin{document}

\title{Accuracy and Performance Evaluation of Quantum, Classical and Hybrid Solvers for the Max-Cut Problem}

\author{Jaka Vodeb}
\email{jaka.vodeb@ijs.si}
\affiliation{Department of Complex Matter, Jožef Stefan Institute, Jamova 39, 1000 Ljubljana, Slovenia}
\affiliation{Department of Physics, Faculty for Mathematics and Physics, University of Ljubljana, Jadranska 19, 1000 Ljubljana, Slovenia}
\orcid{0000-0002-6193-3475}
\author{Vid Er\v{z}en}
\affiliation{Department of Physics, Faculty for Mathematics and Physics, University of Ljubljana, Jadranska 19, 1000 Ljubljana, Slovenia}
\author{Timotej Hrga}
\email{timotej.hrga@aau.at}
\affiliation{Institut für Mathematik, Alpen-Adria-Universität Klagenfurt, Universitätsstraße 65-67, 9020 Klagenfurt am Wörthersee, Austria}
\orcid{0000-0002-4852-1986}
\author{Janez Povh}
\email{janez.povh@rudolfovo.eu}
\affiliation{Rudolfovo - Science and Technology Centre Novo mesto, Podbreznik 15, 8000 Novo mesto, Slovenia}
\affiliation{University of Ljubljana, Faculty of Mechanical Engineering, Aškerčeva 6, 1000 Ljubljana, Slovenia}
\orcid{0000-0002-9856-1476}

\maketitle

\begin{abstract}
  This paper investigates the performance of quantum, classical, and hybrid solvers on the NP-hard Max-Cut and QUBO problems, examining their solution quality relative to the global optima and their computational efficiency. We benchmark the new fast annealing D-Wave quantum processing unit (QPU) and D-Wave Hybrid solver against the state-of-the-art classical simulated annealing algorithm (SA) and Toshiba's simulated bifurcation machine (SBM). Our study leverages three datasets encompassing 139 instances of the Max-Cut problem with sizes ranging from 100 to 10,000 nodes. For instances below 251 nodes, global optima are known and reported, while for larger instances, we utilize the best-known solutions from the literature. 
  Our findings reveal that for the smaller instances where the global optimum is known, the Hybrid solver and SA algorithm consistently achieve the global optimum, outperforming the QPU. For larger instances where global optima are unknown, we observe that the SBM and the slower variant of SA  deliver competitive solution quality, while the Hybrid solver and the faster variant of SA  performed noticeably
worse.  Although computing time varies due to differing underlying hardware, the Hybrid solver and the SBM demonstrate both efficient computation times, while for SA reduction in computation time can be achieved at the expense of solution quality. 
\end{abstract}



\section{Introduction}
\label{intro}

\subsection{Motivation}
We are well into the second quantum revolution, where the knowledge gained from the advent of quantum theory is now being used to build quantum devices. One of the biggest potential benefits from such an undertaking is a new and better way of computing, which has been termed ``quantum computing''. There are two main avenues being developed currently: (i) digital quantum computing, which provides universality in terms of which algorithms can be deployed on such devices, and (ii) quantum annealing (QA), which has only one specific type of algorithm available, but is therefore much more developed in terms of device size and performance. If we only take into account the number of quantum bits or qubits available on the two types of devices, which is by no means a complete performance metric,  the biggest device using digital quantum computing is IBM's Heron with $133$ qubits, while the D-Wave Advantage hosts over $5000$ qubits.

The main driver behind building quantum computers is achieving quantum supremacy, where a quantum computer employing a quantum algorithm significantly outperforms any classical algorithms. This was done by Google in 2019~\cite{arute2019quantum} for a random circuit sampling problem, which holds no practical value. The more realistic goal is a quantum advantage, where a quantum computer outperforms all known classical methods in solving a given problem. 

 According to a recent review~\cite{hauke2020perspectives}, quantum annealers have so far never been able to demonstrate a quantum advantage in practical use cases, which provides motivation for this work. Here, we perform an extensive benchmark test for the Max-Cut and QUBO problems using three standard libraries with problem instances. The result is an assessment of the performance of the Hybrid solver provided by the company D-Wave, as well as their recent more coherent version of quantum processing unit (QPU) involving the fast annealing feature, and their comparison with the performance of (i) the state-of-the-art simulated annealing algorithm for the Max-Cut from~\cite{myklebust2015solving}, and (ii) recent simulated bifurcation machine from~\cite{matsuda2019benchmarking}.
 All these solvers are only approximate, meaning they provide the best solutions they have found without guaranteeing a global optimum. Therefore, we also evaluate their results by comparing them with the best-known values: for problem sizes below 500, these are global optima calculated using exact solvers BiqBin~\cite{gusmeroli2022biqbin} and MADAM~\cite{hrga2021madam}, while for larger problems, we use the best-known values available in the literature.

\subsection{State of the art}

The computational advantage of quantum annealing lies in the ability to traverse the problem's energy landscape via quantum tunneling compared to various classical update schemes used by combinatorial optimization solvers. The landmark study that unveiled the presence of tunneling effects in real-world quantum annealers~\cite{boixo2016computational} was accompanied by an impressive study~\cite{denchev2016computational} by Google, where quantum annealing surpassed simulated annealing by a factor of $\sim10^8$, demonstrating quantum advantage. However, the advantage is limited to specific problem instances, which were designed specifically to boost the performance of quantum annealers in comparison with classical solvers~\cite{denchev2016computational,pang2021potential,tasseff2022emerging,cao2016solving,king2019quantum,mandra2016strengths,koshikawa2021benchmark}. This is not the case in real-world industrial applications, which include problem instances that are even harder to solve than the notorious spin-glass problem~\cite{perdomo2019readiness} with the requirement for embedding an arbitrary problem graph not native to the annealer's quantum processing unit, representing an additional overhead in computational performance~\cite{konz2021embedding}. In principle, quantum annealing provides a theoretical guarantee of finding the ground state due to the quantum adiabatic theorem. However, real devices deviate from this idealized case due to random fields exerted on qubits and random fluctuations of qubit couplers~\cite{zhu2016best}, as well as coupling to an external environment at non-zero temperature~\cite{albash2017temperature}. The effects of both prevent the system from finding the problem solution and even with a perfect device, the annealing time required to find a solution to a generic problem typically scales exponentially with the system size.

Even with a perfect device, the promise of a performance boost from QA comes from the fact that the timescale is set by the quantum instead of classical dynamics. If the smallest energy gap between the ground and first excited state closes polynomially with system size, the time to solution for QA also scales polynomially with system size, which is an exponential advantage over SA. It has been shown that QA outperforms SA~\cite{albash2018demonstration}. However, QA also typically exhibits an exponential closing of the gap, which then provides only an advantage in the prefactor, which can still be significant~\cite{king2021scaling}. It is currently still a matter of debate whether QA exhibits any advantage in relevant industry problems~\cite{hauke2020perspectives}.

Recently, a breakthrough was made by the company D-Wave through faster control of qubits on the quantum processing unit, enabling access to a time scale where qubits are isolated from the environment, thereby avoiding environmental effects \cite{king2022coherent,king2023quantum,king2024computational}. This novel feature remains untested in the context of solving optimization problems.

Recent advancements in quantum annealing hardware have sparked a surge of research and experimentation, driving the quest for enhanced optimization capabilities. Among the forefront contenders in this domain, D-Wave Systems' Advantage quantum annealer stands out, demonstrating notable improvements over its predecessor, the D-Wave 2000Q~\cite{king2019quantum,gonzalez2021garden,dwave2021performanceupdate,mcleod2022benchmarking,willsch2022benchmarking}. Rigorous analyses have unveiled Advantage's prowess in solving a diverse array of optimization problems, such as the stable set problem, constraint satisfaction problems, exact cover problems, maximum cardinality matching problems, random clique problems, 3D lattice problems, and the garden optimization problem, showcasing superior performance metrics such as increased success rates and the ability to tackle larger problem sizes. In particular, Advantage has been shown to triumph over D-Wave 2000Q in scenarios involving complex problems with up to 120 logical qubits, a feat previously beyond the reach of its predecessor~\cite{willsch2022benchmarking}.

In addition to the D-Wave platforms, a plethora of alternative optimization hardware 
has undergone rigorous scrutiny. Fujitsu's digital annealer unit (DAU), virtual MemComputing machines, Toshiba's simulated bifurcation machine (SBM), and various software-based solvers have all been subjected to comprehensive scaling analyses and prefactor evaluations. These investigations show that quantum annealing exhibits the worst time-to-solution scaling in comparison to other platforms~\cite{kowalsky20223}. In particular, quantum annealing hardware is not well-suited for finding degenerate ground-state configurations associated with a particular instance~\cite{mandra2017exponentially,zhu2019fair}. A comparison with measurement-feedback coherent Ising machines (CIMs) based on optical parametric oscillators on two problem classes, the Sherrington-Kirkpatrick (SK) model and Max-Cut showed that the D-Wave quantum annealer outperforms CIMs on Max-Cut on cubic graphs~\cite{hamerly2019experimental}. On denser problems, however, there is a several order of magnitude time-to-solution difference for instances with over 50 vertices~\cite{hamerly2019experimental}. Therefore, CIMs are not a clear replacement for quantum annealers, and they also suffer from underperformance in comparison to classical combinatorial optimization solvers~\cite{sharma2023augmenting}.

The intrinsic heuristic nature of quantum annealers prevents knowing whether the best solution produced by the device is also the global minimum, which can be used to assess the accuracy of the annealer or, in other words, how close the annealer comes to actually solving the proposed problem. Exact solvers have already been used to construct exhaustive libraries of problems and the best-found solutions. These instances are readily available for benchmarking quantum annealers and have not been used so far to the best of our knowledge. A different approach would be to construct a problem specifically tailored to the quantum chip topology where the solution is known by design~\cite{pelofske2024increasing}, but this is not the focus of our work.

In general, better quality qubits, higher connectivity, and more tunability is needed to answer the question if quantum annealing will ever truly outperform specialized silicon technology combined with efficient heuristics for optimization and sampling applications~\cite{katzgraber2018viewing}, thereby making the recent hardware update enabling faster annealing times an important topic of benchmarking. On the other hand, developing software methods that make the best use of the quantum technology available today~\cite{montanez2023improving} is crucial in pushing the boundaries of quantum annealing hardware capabilities, given the great potential it holds in the field of machine learning~\cite{radovic2018machine,zlokapa2020quantum}.

\subsection{Our contribution}
The studies referenced above, which examine the performance of quantum annealers on combinatorial optimization problems, focus primarily on questions concerning the largest instance sizes that a given solver can address, the computation time required, and the quality of solutions compared to those obtained with other solvers. Notably, most combinatorial optimization problems explored in the literature are NP-hard, indicating that if we want to solve any instance of these problems, i.e., find a global optimum, we will face substantial theoretical and practical complexity. If near-optimal solutions are sufficient, practical heuristic algorithms can often provide them efficiently.

This naturally raises the question of how close the solutions produced by various solvers are to the global optimum. Therefore, the main objectives of this paper are: (i) to assess the proximity of solutions obtained for the Max-Cut and QUBO problems by the D-Wave fast annealing quantum processing unit (QPU) and the D-Wave Hybrid solver to the global optimum, and (ii) to benchmark the performance of these solvers against the state-of-the-art simulated annealing algorithm for the Max-Cut from  ~\cite{myklebust2015solving} and the recently introduced Simulated Bifurcation Machine (SBM) from~\cite{matsuda2019benchmarking}.

\vspace{1em}
The main contributions of this paper are:
\begin{itemize}
    \item We prepared three benchmark datasets for the Max-Cut and QUBO problems, encompassing 139 instances with sizes ranging from 100 to 10,000. For the instances with sizes below 251 (this includes the first dataset and half of the second dataset), global optima are provided by using exact solvers BiqBin ~\cite{gusmeroli2022biqbin} and MADAM~\cite{hrga2021madam}, while for larger instances, we include the best-known solutions from the literature. These datasets are publicly available as open data in our GitHub repository \cite{erzen_2024_14290290} and are formatted for seamless use in further benchmarking studies.
    \item We (approximately) solved the instances from the benchmark datasets by the D-Wave fast annealing quantum processing unit (QPU) and D-Wave hybrid solver (Hybrid).
    \item We (approximately) solved the same instances by the state-of-the-art simulated annealing solver for Max-Cut from~\cite{myklebust2015solving}. We used two variants, SA1 and SA2, of this algorithm that utilize different annealing schedules, i.e., at which initial temperature we start and how big step sizes we make till cooldown. The variant SA1 is more aggressive with the step size, while SA2 spends more time in the search space. 
    \item We conduct an in-depth comparison of solvers in terms of the proximity of their solutions to the global optima. On the first dataset, which contains small instances, we compared QPU, Hybrid,  SA1, and SA2 solvers in terms of proximity to optimum solutions. On the second data set, we compared only the Hybrid solver against the simulated annealing solvers, since the instances exceeded the capacity of the QPU. We compared the proximity of computed solutions to optimum or best-known solutions. On the third dataset, which contains the largest instances, we compared the  Hybrid solver against the Toshiba Simulated Bifurcation Machine and the simulated annealing solvers SA1 and SA2. We considered the quality of solutions and the computational times required for solution computation.
\end{itemize}

 Our analysis shows that on the first two datasets, the Hybrid solver and both variants Simulated Annealing solvers compute global optima whenever they are known or compute solutions with objective values equal to the best-known values from the literature, while QPUs solutions are usually quite far from the optima or best-known solutions.  
 
 On the third dataset, for which the global optimum is not known, we observe that the SBM and simulated annealing solver SA2 are very competitive in terms of the quality of solutions, while the Hybrid solver and the SA1 solver compute   
 solutions that are further away from the optimum or the best known solutions.
 The corresponding computing times, which are actually hard to compare due to different underlying hardware, indicate that Hybrid,  SBM are also very fast, while for simulated annealing, we must find a balance between the quality and the computing time since more computing time gives much better solutions, as is the case with SA2.

\section{Max-Cut and QUBO problems}
In this paper, we will consider only undirected weighted graphs without loops and multi-edges, i.e., simple weighted graphs. Each graph $G$ is defined by a set of vertices $V=\{v_1,v_2,\ldots,v_n\}$, set of edges $E=\{\{v_i,v_j\}\mid v_i\neq v_j ~\mbox{and}~  v_i~\mbox{is connected with}~v_j\}$ and the matrix of edge weights $A$, where   $a_{ij}$ is the weight on the edge between $v_i$ and $v_j$, if $v_i$ and $v_j$ are  connected, otherwise we have $a_{ij}=0.$ If the graph is unweighted, we simply have $a_{ij}\in \{0,1\}.$ We usually represent a simple (weighted or unweighted)  graph $G$ by the matrix $A$, which is called the adjacency matrix of $G$.

Given an edge-weighted undirected graph, the Max-Cut problem asks to find the size of the largest cut in a given graph. A cut is a set of edges
that separates the vertices $V$ into two disjoint sets $V_1$ and $V_2$,
such that $V_1 \subseteq V$ and $V_2 = V \backslash V_1$, and the cost of a cut is
defined as the sum of all weights $a_{ij}$ of edges connecting vertices in $V_1$ with vertices in $V_2$.
The Max-Cut problem is a fundamental and well-known example of NP-hard combinatorial optimization problems~\cite{karp2010reducibility}, which are notoriously hard to solve. It has a surprisingly large number of important practical applications, for example, in VLSI design~\cite{barahona1988application} and in the theory of spin glasses in physics~\cite{grotschel1987calculating, liers2004contributions}. 

The cost of a cut can be computed by assigning the binary values to express
which subset the vertex $v_i$ belongs to $x_i \in \{-1,1\}$, i.e.,
every partition of $V$ can be represented by a vector $\x\in\{-1,1\}^n$ where $x_i=-1$ if and only if the $i$-th vertex $v_i$ belongs to $V_1$. 
Let $\mbox{cut}(V_1,V_2)$ denote the cost of a cut. With this notation, we have
 
  \[ \mbox{cut}(V_1,V_2)=\sum_{i \in V_1, j \in V_2} a_{ij} = \sum_{i < j}a_{ij}\frac{1-x_ix_j}{2} = \frac{1}{4}\x^TL\x,\]
  where $L$ is the Laplacian matrix of graph $G$, defined by  $L=\Diag(A\e)-A$. We denoted by $\e$ the vector of all ones and by $\Diag(A\e)$ the diagonal matrix,  having on the main diagonal the vector $A\e$.
  
Therefore, the Max-Cut problem can be formulated as a binary quadratic program:
\begin{align}
  \label{eqn:MC}  \tag{MC}
  \begin{split}
  \max~~ & \frac{1}{4}\x^T L \x\\
  \st~~  &\x \in \{\pm 1\}^{n}.
  \end{split}
\end{align}

Many heuristics have been proposed for Max-Cut. See Dunning et~al.~\cite{dunning2018works} for a recent overview. There have also been various articles on exact solution approaches, see e.g. Jünger et~al.~\cite{junger2021quantum} for an up-to-date overview. 

Goemans and Williamson~\cite{goemans1995improved} devised a polynomial-time randomized rounding algorithm based on a semidefinite relaxation that computes a solution whose objective value is within a factor $\alpha \approx 0.87856$ of the optimum whenever the weights of the underlying graph are nonnegative. Khot et~al.~\cite{khot2007optimal} proved that, if the Unique Games Conjecture is true, then there is no polynomial-time algorithm that always computes a solution
whose objective value is within a constant factor of optimum, where this factor is larger than $\alpha$.

The Max-Cut problem can be considered a special instance of the  classical Ising problem of ferromagnetism in statistical mechanics, where one tries to minimize the classical Hamiltonian of system of spins:
\begin{equation}\label{eq:Ham1}\tag{Hamiltonian}
H = \frac{1}{2}\sum_{i,j}J_{ij}\sigma_i\sigma_j + \sum_ih_i\sigma_i
\end{equation}
where $\sigma_i = \pm 1$ represent the values of the Ising spins, $J_{ij}$ are the entries of the spin-spin coupling matrix, and the values $h_i$ represent the influence of the external magnetic field on the spins.
The minimization is performed over all possible spin configurations $\sigma$.
It is easy to prove that the Ising problem and the Max-Cut problem are actually equivalent. The Max-Cut problem is an instance of the Ising problem without an external magnetic field ($h_i=0,\forall i$), and the Ising problem becomes the Max-Cut problem if we augment the matrix $J$ by adding the vector $\frac{1}{2}h$ to the first row and first column and increase the dimension of the problem by one.

We also introduce another quadratic optimization problem, closely related 
to the Max-Cut problem:
The Quadratic Unconstrained Binary Optimization (QUBO) problem. It  can be formulated as follows: for a given symmetric $n\times n$ matrix $Q$ we want to solve 
\begin{align}
  \label{eqn:QUBO}  \tag{QUBO}
  \begin{split}
  \min~~ & \x^T Q \x\\
  \st~~  &\x \in \{0,1\}^{n}.
  \end{split}
\end{align}
It is not difficult to show that solving \eqref{eqn:QUBO} is equivalent to solving \eqref{eqn:MC}. QUBO and Max-Cut are equivalent problems, as there exists a linear transformation between them. In particular, any QUBO problem on $n$ variables can be transformed into an equivalent Max-Cut instance on $n + 1$,  and any Max-Cut instance on a graph $(V, E)$ can be formulated
as a QUBO instance with $n = |V| - 1$, see e.g.~\cite{barahona1989experiments}.
The reason for the prominence of QUBO formulation in quantum computing is that many current quantum hardware platforms have a natural connection to QUBO. 
In Section~\ref{subsec:datasets} we describe four datasets that are used for benchmarking different solvers. Three of them have QUBO origins, i.e. they were originally formulated as QUBO and later reformulated as Max-Cut. In this paper, we always use the Max-Cut formulations.

\section{Description of used solvers for the Max-Cut problem}

The instances of the Max-Cut problem that we consider in this paper are well-known benchmark instances, extensively studied in the literature. 
For many of them, the optimum solution is known, mostly computed by some variant of complete enumeration algorithms, like the Branch-and-bound or Branch-and-cut algorithms. For all instances with $n\le 251$ we computed the optimum value by our solvers BiqBin~\cite{gusmeroli2022biqbin} and MADAM~\cite{hrga2021madam}. We refer in the tables below only to BiqBin solver.
For instances with $250 < n \le 500$, our solvers gave us only good solutions (often equal to the best known solutions), without a certificate that these are optimum solutions.

Besides the exact solvers, we also applied to the benchmark instances four approximate solvers.  The D-Wave fast annealing quantum processing unit solver (QPU) and the out-of-the-box Hybrid solver (Hybrid) are two approximate solvers provided by D-Wave Systems, which are accessible as cloud solutions. 
The QPU solves a problem instance in an analogous way to simulated annealing, where temperature fluctuations are replaced with quantum fluctuations. 
  The classical Hamiltonian \eqref{eq:Ham1}, representing the Max-Cut, is encoded into the extended Hamiltonian
\begin{align}
  \label{eqn:DWH}  \tag{D-Wave Hamiltonian}
  \begin{split}
 H = -\frac{1}{2}A(s)\sum_{i}\sigma_i^x + \frac{1}{2}B(s)\left(\sum_{i<j}J_{i,j}\sigma_i^z\sigma_j^z + \sum_ih_i\sigma_i^z\right),
  \end{split}
\end{align}
where  $\sigma_i^z$ and $\sigma_i^x$ are Pauli operators, which describe the dynamics of the qubit representing $x_i$. The parameter $s\in[0,1]$ is used to drive the values of scalar functions $A(s)$ and $B(s)$ within a time frame $t_a$. After the Hamiltonian is encoded on the quantum annealer, $s$ is set to $0$, where $A(0)$ is non-zero and $B(0)=0$. Then within time $t_a$ the value of $s$ goes to $1$, where $A(1)=0$ and $B(1)$ is non-zero. This procedure is analogous to the simulated annealing (SA) method, where the temperature is first set to infinity (or a high enough value above the phase transition in practice) and is then reduced to $0$. 

According to the publicly available description of the Hybrid\footnote{https://cloud.dwavesys.com/leap/example-details/285687969/}, it is comprised of an iterated loop consisting of several classical optimization methods as well as the QPU. The classical optimization methods include an interruptable simulated annealing sampler, tabu sampler and a parallel tempering sampler. These methods are run simultaneously in an unspecified way within the Python method {\tt LeapHybridSampler()}, which we employed for this work. However, the providers claim that this is a method optimized for a wide variety of optimization problems, not for any specific applications, and encourage users to develop their own workflows in order to optimize the method for their own specific purposes. We stuck to the out-of-the-box method due to the variety of our input problem instances.

The sub-method of the Hybrid, which employs the QPU, is comprised of an energy decomposer that reduces the input graph to smaller graphs, solves them, and then puts them back together. The public does not have access to the details of the exact methods used for the decomposition and reconstruction of the input graph. When the QPU returns a solution to the whole input graph, the whole iteration loop terminates.
With the QPU we could solve only instances with $n\le 151$, while Hybrid could solve all instances.

Additionally, we also used classical simulated annealing (SA) which was independently developed by Kirkpatrick et al.~\cite{kirkpatrick1983optimization} and Čern\'y\cite{vcerny1985thermodynamical}.
Simulated annealing works by taking a random walk in the space of feasible solutions. It considers a random local modification of the current solution. 
If that modification leads to a better solution, it is always accepted. If not, its new modification is made with some probability that decreases exponentially with the amount of objective degradation. By accepting points that increase the objective value, the algorithm avoids being trapped in local minima, and is able to explore globally far more possible solutions.

The implementation of the simulated annealing for solving Max-Cut problems that we used for computations is from~\cite{myklebust2015solving}. 
The reason for this is that in their version, each iteration is computationally efficient, and this allows the code to explore a substantially large number of feasible solutions. Furthermore, the parameters that determine the initial temperature and decrement heavily influence the search time and have an impact on the runtime of the algorithm. 
We used two variants that employ different annealing schedules. For the SA1 version, the initial temperature was set to 10,000 with the decrement step of $2\mathrm{e}{-4}$ in order to get the runtimes that are comparable with D-Wave. In~\cite{myklebust2015solving}, the author reports that by changing the initial temperature to 40,000 and the decrement to heat to $5\mathrm{e}{-6}$, the SA finds a feasible solution to instance G35 with objective value 7685, which is better than those found by any other heuristic. Motivated by this and to get better feasible solutions, we set the starting value of the temperature to 40,000 with the decrement step of $2\mathrm{e}{-6}$ for version SA2.

\section{Numerical results}
In this section, we introduce three benchmark datasets and present numerical results obtained from all solvers capable of providing solutions. For instances for $n \leq 151$, we employed the exact solvers BiqBin and MADAM, the fast annealing quantum solver QPU, the hybrid solver, and two simulated annealing solvers, SA1 and SA2. For instances with $n = 251$, the QPU solver exceeded its capacity, and for instances with $n > 251$, even the exact solvers could no longer handle the computations.

\subsection{Datasets}\label{subsec:datasets}

We  performed computations  on three Max-Cut datasets, where two of them were originally formulated as QUBO:
\begin{itemize}
    
    \item The  instances of  Billionnet and Elloumi  from  \cite{wiegele2007biq,biqmaclib} with $n = 100,120,150$, denoted by {\tt be100.1,\ldots,be100.10, be120.3.1,\ldots,be120.3.10,\\ be120.8.1,\ldots,be120.8.10, be150.3.1,\ldots,be150.3.10}, and \\
    {\tt be150.8.1,\ldots,be150.8.10}. These instances were originally formulated as QUBO, but we used the Max-Cut reformulations. The source files for our computations were   {\tt be100.1.sparse.mc,\ldots,be150.8.10.sparse.mc},\\ computed from the data files 
    {\tt be100.1.sparse,\ldots,be150.8.10.sparse} with a code {\tt qp2mc}, written by prof.~Angelika Wiegele. The files {\tt *.sparse} were taken from the  BiqMac site\footnote{\url{https://biqmac.aau.at/biqmaclib.html}}. The new files {\tt *.sparse.mc} are available on our GitHub repository~\cite{erzen_2024_14290290}. 
    
    We created matrices $\frac{1}{4}L$ from these files and named them as 
    {\tt be100.1\_L.txt, \ldots,} \\ {\tt  be150.8.10\_L.txt}. From these matrices, we extracted Python dictionaries needed as input files for the D-Wave solvers. From each matrix $\frac{1}{4}L$, we extracted the lower triangular matrix (including the diagonal) and multiplied it with $(-1)$ since the D-Wave solvers solve minimization problems. The resulting matrices are available at the GitHub repository  \cite{erzen_2024_14290290} as Python dictionaries with names \\
    {\tt be100.1\_py.txt,\ldots,be150.8.10\_py.txt.}
    Additionally,  for each instance  we also provided an optimum cut vector \\
    {\tt be100.1\_opt\_cut.txt,\ldots,be150.8.10\_opt\_cut.txt}  and the corresponding optimum values   {\tt be100.1\_opt\_value.txt,\ldots, 
 }  {\tt be150.8.10\_opt\_value.txt}. 
    Note that the numbers in the names of the instances tell the size of the original QUBO formulation, but after the reformulation into the Max-Cut problem, the size of the instance is increased by 1, so, e.g., the Max-Cut formulation of instance {\tt be100.1} has $n=101$.
    
   \item Instances from the Beasley collection \cite{wiegele2007biq,biqmaclib} with $n = 250$ and $n=500$, denoted by {\tt bqp250.1,\ldots,bqp250.10} and {\tt bqp-500.1,\ldots,bqp500.10.} These instances were also originally created as QUBO problems, but we use the Max-Cut reformulation.
   Like with the {\tt be} instances, the input files were the files {\tt bqp250-1.sparse.mc,\ldots, bqp500-10.sparse.mc.} Again, we created from  them the matrices $\frac{1}{4}L$ and the Python dictionaries with lower negative triangular part of  $\frac{1}{4}L$ (including diagonals). 
  %  \item \textcolor{red}{Za te primere nimamo vključenih numeričmnih rezultatov. Manjkajo SA1 in SA2.}
   % Instances of QUBO, denoted by {\tt gka1e, gka1f, gka2e, gka2f, \ldots gka5e, gka5f} from \cite{glover1998adaptive}. These instances were actually introduced already by  Pardalos and Rodgers in \cite{pardalos1990computational}. We denoted them by {\tt gka}, 
   % since they were selected and described in details in    \cite{glover1998adaptive}.
   % We use the Max-Cut formulations, so he input matrices for these instances were matrices {\tt gka1e.mc,\ldots,gka5f.mc}. Again, we created out of them the matrices $\frac{1}{4}L$ and the corresponding Python dictionaries, containing negative lower triangular parts, including the diagonals. Resulting matrices are named {\tt gka1e\_L.txt,\ldots,gka5f\_L.txt} and {\tt gka1e\_py.txt,\ldots,gka5f\_py.txt}.
   %\end{comment}
 As with {\tt be} instances, the dimension of the Max-Cut formulations of {\tt bqp} instances is 1 greater than that of the QUBO formulation.
 As with the {\tt be} instances, we provide in the GitHub repository \cite{erzen_2024_14290290} for each instance the input files, an optimal cut vector and the optimum value when $n=251$, and a cut vector giving the best-known cut, including the value of the best-known cut when $n=501$.
   
\item The G-dataset from Yinyu Ye\footnote{http://web.stanford.edu/~yyye/yyye/Gset/}, as used in \cite{matsuda2019benchmarking}. 
    The set of Max-Cut instances was generated by Helmberg and Rendl \cite{helmberg2000spectral} using G. Rinaldi’s graph generator \verb|rudy|.
    It consists of 69 matrices labeled as G1, G2,\ldots, and G72 (G68, G69, and G71  are missing for a reason not known to us). 
Here, our workflow was slightly different. We took as input the edge list format of graph adjacency matrices {\tt G1.txt,\ldots,G72.txt} from Yinyu Ye's web page and created out of them the matrices $L$ (not $\frac{1}{4}L$) and the corresponding Python dictionaries with negative lower triangular parts, without the diagonal.
Results were stored with names {\tt G1\_L.txt,\ldots,G72\_L.txt} and {\tt G1\_py.txt,\ldots,G72\_py.txt}. In the GitHub repository \cite{erzen_2024_14290290}, we provide for each instance the original data file,  the $L$ matrix, the Python dictionary, a cut vector yielding the best-known cut value, and the best-known cut value. 
\end{itemize}

%Beside the  matrices $L$ or $\frac{1}{4}L$, as described above, and the associated Python dictionaries, we publish for each instance on our GitHub repository also  the cut vector $ \x$, yielding the best computed value, and the best computed value.

\subsection{Numerical results on benchmark datasets}

Our numerical results were obtained using a variety of computing infrastructures. All exact values, computed by BiqBin or MADAM, were computed on the ROME partition of the HPC system at the University of Ljubljana, Faculty of
Mechanical Engineering. We used a serial version of the solvers, which ran on
single computing nodes, each consisting of 2 AMD EPYC 7402 24-core processor, 128 GB DDR4-3200 RAM, and 1 TB NVMe SSD. The tests with classical simulated annealing were performed on a machine with 2.8 GHz Dual-Core Intel Core i7 with 8GB of RAM.

On the {\tt be} instances, which are small enough, we executed the exact solver BiqBin, the fast annealing quantum processing unit solver (QPU), the out-of-the-box Hybrid solver (Hybrid), and the simulated annealing solvers SA1 and SA2. The results are given in Table~\ref{tab:res1}. The first column of each row gives the name of the
underlying graph, the second column contains the size of the problem, and the third column, labeled BiqBin, gives the optimum solutions of the Max-Cut problem obtained by BiqBin or MADAM solvers (both were used for these computations). In the next four columns, we report solutions provided by QPU, Hybrid,  SA1, and SA2 solvers. We can see that Hybrid, SA1, and SA2 returned the optimum values for all test instances.

Table~\ref{tab:res2} contains numerical results for larger instances. For instances with $n=251$, the optimal values for Max-Cut calculated by BiqBin and the values provided by Hybrid,  SA1, and SA2 are given. For $n=501$, the BiqBin column contains the best values calculated by BiqBin, without a certificate that it has found the optimum. Like in Table~\ref{tab:res1}, for $n=251$, Hybrid,  SA1, and SA2 all find the optimal solutions. For $n=501$, all four solvers find equivalent solutions, i.e. the cuts they found have the same values.

In Tables~\ref{tab:res3} and~\ref{tab:res4} we report the numerical results for the G data set. These instances are too large to use BiqBin and QPU, so we applied only the Hybrid solver and both variants of simulated annealing solvers, SA1 and SA2. Additionally, we report results from Toshiba Simulated  Bifurcation Machine
solver (SBM), taken from~\cite{matsuda2019benchmarking}.
In addition to the computed value of the Max-Cut problem, we also give the computational times in seconds.
We are aware that the times are not directly comparable, as the infrastructures on which we performed the calculations were very different. Nevertheless, we have included them to show that we can obtain very good solutions in a relatively short time. Since we do not know the optimal solutions, we only highlight the best-computed solution among all four solvers, which are also the best-known solutions so far. We can see that SBM and SA2 are very competitive in terms of the quality of the solutions. However, based on the reported computation times ~\cite{matsuda2019benchmarking}, SA2 takes much more time (factor 10000).
The Hybrid solver is  worse than SA2 in terms of computed values, but computational times are much closer to those of SBM.
The solver SA1 takes about the same time as Hybrid and often slightly more than SBM, but computed values are, in most cases, worse than with the other three solvers. 
%We colored red two instances where SA1 computed better solutions than SA2.

A summary of all results is given in Table~\ref{tab:res5}, from which it is clear that datasets \texttt{be} and \texttt{bqp} were not demanding enough to differentiate between the Hybrid, SA1, and SA2  solvers - all of them produced solutions of the same quality for all these problems. Dataset \texttt{G}, on the other hand, was able to distinguish the Hybrid solver, SBM and both SA algorithms according to the quality of solutions. SBM produced the greatest number of best solutions and the greatest number of solutions that were strictly better compared to the other solvers.  The solver SA2 was close SBM, while the Hybrid and SA1 solvers performed noticeably worse. SA1 exhibited the worst performance among all.
It is interesting that SA2, which is a traditional simulated annealing algorithm, adapted to Max-Cut, outperforms the very advanced black-box Hybrid solver on the \texttt{G} dataset.

\begin{widetext}
\begin{table}[H]\small 
\centering
\caption{A comparison between the optimal values obtained with the BiqBin exact solver, and the best objective values obtained with D-Wave's fast annealing QPU and the D-Wave Hybrid solver, and the Simulated Annealing heuristics SA1 and SA2.}\label{tab:res1}
\begin{tabular}{lcccccc}
\toprule
Instance & $n$ & BiqBin & QPU & Hybrid & SA1 & SA2 \\
\midrule
be100.1 & 101 & -19412 &-19258 &-19412 & -19412 & -19412 \\ 
be100.2 & 101 & -17290 &-17174 &-17290 & -17290 & -17290\\
be100.3 & 101 & -17565 &-17389 &-17565 & -17565 & -17565 \\
be100.4 & 101 & -19125 &-18988 &-19125 & -19125 & -19125 \\
be100.5 & 101 & -15868 &-15661 &-15868 & -15868 & -15868\\
be100.6 & 101 & -17368 &-17368 &-17368 & -17368 & -17368 \\
be100.7 & 101 & -18629 &-18368 &-18629 & -18629 & -18629 \\
be100.8 & 101 & -18649 &-18641 &-18649 & -18649 & -18649 \\
be100.9 & 101 & -13294 &-12807 &-13294 & -13294 & -13294 \\
be100.10 & 101& -15352 &-15154 &-15352 & -15352 & -15352 \\
\midrule
be120.3.1 & 121 & -13067 &-12342 &-13067 & -13067 & -13067\\
be120.3.2 & 121 & -13046 &-12726 &-13046 & -13046 & -13046 \\
be120.3.3 & 121 & -12418 &-12156 &-12418 & -12418 & -12418 \\
be120.3.4 & 121 & -13867 &-13381 &-13867 & -13867 & -13867 \\
be120.3.5 & 121 & -11403 &-11237 &-11403 & -11403 & -11403 \\
be120.3.6 & 121 & -12915 &-12582 &-12915 & -12915 & -12915\\
be120.3.7 & 121 & -14068 &-13946 &-14068 & -14068 & -14068 \\
be120.3.8 & 121 & -14701 &-14016 &-14701 & -14701 & -14701 \\
be120.3.9 & 121 & -10458 &-10230 &-10458 & -10458 & -10458 \\
be120.3.10& 121 & -12201 &-11375 &-12201 & -12201 & -12201 \\
\midrule
be120.8.1 & 121 & -18691 &-18372 &-18691 & -18691  & -18691\\
be120.8.2 & 121 & -18827 &-18161 &-18827 & -18827 & -18827 \\
be120.8.3 & 121 & -19302 &-18561 &-19302 & -19302 & -19302\\
be120.8.4 & 121 & -20765 &-19536 &-20765 & -20765 & -20765\\
be120.8.5 & 121 & -20417 &-20224 &-20417 & -20417 & -20417 \\
be120.8.6 & 121 & -18482 &-18132 &-18482 & -18482 & -18482 \\
be120.8.7 & 121 & -22194 &-21431 &-22194 & -22194 & -22194\\
be120.8.8 & 121 & -19534 &-17832 &-19534 & -19534 & -19534\\
be120.8.9 & 121 & -18195 &-17511 &-18195 & -18195 & -18195 \\
be120.8.10 &121 & -19049 &-18604 &-19049 & -19049 & -19049 \\
\midrule
be150.3.1  &151 & -18889 &  -17804 &-18889 & -18889 & -18889\\
be150.3.2  &151 & -17816 &  -17103 &-17816 & -17816 & -17816\\ 
be150.3.3  &151 & -17314 &  -16609 &-17314 & -17314 & -17314 \\ 
be150.3.4  &151 & -19884 &  -18231 &-19884 & -19884 & -19884  \\ 
be150.3.5  &151 & -16817 &  -15904 &-16817 & -16817 & -16817 \\ 
be150.3.6  &151 & -16780 &  -16346 &-16780 & -16780 & -16780 \\
be150.3.7  &151 & -18001 &  -16156 &-18001 & -18001 & -18001 \\
be150.3.8  &151 & -18303 &  -17054 &-18303 & -18303 & -18303 \\
be150.3.9  &151 & -12838 &  -11728 &-12838 & -12838 & -12838 \\
be150.3.10 &151 & -17963 &  -16117 &-17963 & -17963 & -17963 \\
\midrule
be150.8.1  &151 & -27089 &  -24231 &-27089 & -27089 & -27089 \\
be150.8.2  &151 & -26779 &  -25912 &-26779 & -26779 &-26779 \\
be150.8.3  &151 & -29438 &  -26448 &-29438 & -29438 &-29438 \\
be150.8.4  &151 & -26911 &  -24586 &-26911 & -26911 &-26911 \\
be150.8.5  &151 & -28017 &  -26874 &-28017 & -28017 &-28017 \\
be150.8.6  &151 & -29221 &  -27481 &-29221 & -29221 &-29221 \\
be150.8.7  &151 & -31209 &  -28275 &-31209 & -31209 &-31209 \\
be150.8.8  &151 & -29730 &  -27397 &-29730 & -29730 &-29730 \\
be150.8.9  &151 & -25388 &  -22949 &-25388 & -25388 &-25388 \\
be150.8.10 &151 & -28374 &  -27251 &-28374 & -28374  &-28374
\\ \bottomrule
\end{tabular}
\end{table}
\end{widetext}

\begin{table}[H]\small 
\centering
\caption{A comparison between the best objective values obtained with the BiqBin solver, the D-Wave Hybrid solver, and the Simulated Annealing heuristics SA1 and SA2.}\label{tab:res2}
\begin{tabular}{lccccc}
\toprule
Instance & $n$ & BiqBin & Hybrid & SA1 & SA2 \\
\midrule
bqp250-1 & 251 & -45607 &-45607 & -45607 & -45607\\
bqp250-2 & 251 & -44810 &-44810 & -44810 & -44810 \\
bqp250-3 & 251 & -49037 &-49037 & -49037 & -49037\\
bqp250-4 & 251 & -41274 &-41274 & -41274 & -41274 \\
bqp250-5 & 251 & -47961 &-47961 & -47961 & -47961 \\
bqp250-6 & 251 & -41014 &-41014 & -41014 & -41014 \\
bqp250-7 & 251 & -46757 &-46757 & -46757 & -46757 \\
bqp250-8 & 251 & -35726 &-35726 & -35726 & -35726\\
bqp250-9 & 251 & -48916 &-48916 & -48916 & -48916 \\
bqp250-10 &251 & -40442 &-40442 & -40442 & -40442\\
\midrule
bqp500-1 & 501 & -116586 &-116586 &  -116586 & -116586 \\
bqp500-2 & 501 & -128339 &-128339 &  -128339 & -128339\\
bqp500-3 & 501 & -130812 &-130812 &  -130812 &  -130812 \\
bqp500-4 & 501 & -130097 &-130097 &  -130097 &  -130097\\
bqp500-5 & 501 & -125487 &-125487 &  -125487 &  -125487 \\
bqp500-6 & 501 & -121772 &-121772 &  -121772 &  -121772 \\
bqp500-7 & 501 & -122201 &-122201 &  -122201  &  -122201 \\
bqp500-8 & 501 & -123559 &-123559 &  -123559 &  -123559 \\
bqp500-9 & 501 & -120798 &-120798 &  -120798 &  -120798\\
bqp500-10 & 501 & -130619&-130619 &  -130619 &  -130619\\
\bottomrule
\end{tabular}
\end{table}

\vspace{2em}

\begin{table}[H]
\centering
\caption{A comparison between the best solutions obtained with the Simulated Annealing heuristic (SA), Toshiba Simulated  Bifurcation Machine solver (SBM), and the D-Wave Hybrid solver.}\label{tab:res3}
\begin{adjustbox}{width=1\textwidth}
\begin{tabular}{cccccccccc}
\toprule
Instance & $n$ & Hybrid & $t_{Hybrid} (s)$ & SBM & $t_{SBAM} (s)$  & SA1 & $t_{SA1} (s)$ & SA2 & $t_{SA2} (s)$ \\
\midrule
G1  & 800  & {\textbf{11624}}  &  3.00  & {\textbf{11624}}  & 0.02  & 11615 & 2.66 & {\textbf{11624}} & 672.70  \\
G2  & 800  & {\textbf{11620}}  &  3.00  & {\textbf{11620}}  & 0.13  & 11604 & 2.64 & {\textbf{11620}} & 672.01  \\
G3  & 800  & {\textbf{11622}}  &  2.99  & {\textbf{11622}}  & 0.09  & 11610 & 2.82 & {\textbf{11622}} & 675.80  \\ 
G4  & 800  & {\textbf{11646}}  &  3.00  & {\textbf{11646}}  & 0.02  & 11626 & 2.75 & {\textbf{11646}} & 671.84  \\ 
G5  & 800  & {\textbf{11631}}  &  3.00  & {\textbf{11631}}  & 0.58  & 11618 & 2.64 & {\textbf{11631}} & 674.53  \\
G6  & 800  & {\textbf{2178}}   &  3.00  & {\textbf{2178}}   & 0.02  & {\textbf{2178}}  & 1.64 & {\textbf{2178}}  & 556.42  \\ 
G7  & 800  & {\textbf{2006}}   &  3.00  & {\textbf{2006}}   & 0.02  & {\textbf{2006}}  & 1.62 & {\textbf{2006}}  & 569.27  \\ 
G8  & 800  & {\textbf{2005}}   &  3.00  & {\textbf{2005}}   & 0.02  & 1999  & 1.62 & {\textbf{2005}}  & 572.57 \\
G9  & 800  & 2053   &  3.00  & {\textbf{2054}}   & 0.02  & 2046  & 1.57 & {\textbf{2054}}  & 567.68 \\
G10 & 800  & {\textbf{2000}}   &  3.00  & {\textbf{2000}}   & 0.06  & 1998  & 2.16 & 1999  & 572.10\\ 
G11 & 800  & 562    &  3.00  & {\textbf{564}}    & 0.01  & 558   & 1.49 & {\textbf{564}}   & 571.73\\
G12 & 800  & 554    &  3.00  & {\textbf{556}}    & 0.02  & 552   & 1.48 & {\textbf{556}}   & 573.89\\
G13 & 800  & 580    &  2.99  & {\textbf{582}}    & 0.02  & 576   & 1.43 & {\textbf{582}}   & 571.59\\ 
G14 & 800  & 3058   &  3.00  & {\textbf{3063}}   & 10.00 & 3059  & 1.86 & 3062  & 611.48\\ 
G15 & 800  & 3045   &  3.00  & {\textbf{3050}}   & 0.06  & 3035  & 1.87 & 3049  & 615.67\\ 
G16 & 800  & 3047   &  2.99  & {\textbf{3052}}   & 0.18  & {\textbf{3052}}  & 1.86 & {\textbf{3052}}  & 603.63 \\ 
G17 & 800  & 3043   &  3.00  & {\textbf{3047}}   & 4.64  & 3042  & 1.80 & 3045  & 606.08\\ 
G18 & 800  & 988    &  2.98  & {\textbf{992}}    & 0.16  & 986   & 1.54 & 988   & 575.92\\ 
G19 & 800  & {\textbf{906}}    &  2.99  & {\textbf{906}}    & 0.02  & {\textbf{906}}  & 1.52 &  903    & 579.26\\ 
G20 & 800  & {\textbf{941}}    &  3.00  & {\textbf{941}}    & 0.01  & {\textbf{941}}   & 1.55 & {\textbf{941}}   & 577.41\\ 
G21 & 800  & 928    &  3.01  & {\textbf{931}}    & 0.08  & 927   & 1.54 & 927   & 578.4\\ 
G22 & 800  & 13351  &  5.22  & {\textbf{13359}}  & 0.20  & 13099 & 2.71 & {\textbf{13359}} & 707.93\\ 
G23 & 800  & 13327  &  5.22  & 13342  & 9.99  & 13077 & 2.70 & {\textbf{13344}} & 708.00\\ 
G24 & 2000 & 13330  &  5.22  & {\textbf{13337}}  & 0.07  & 13081 & 2.73 & {\textbf{13337}} & 706.96\\ 
G25 & 2000 & 13334  &  5.22  & {\textbf{13340}}  & 0.45  & 13088 & 2.72 & 13333 & 707.71\\ 
G26 & 2000 & 13322  &  5.22  & {\textbf{13328}}  & 0.07  & 13077 & 2.73 & 13323 & 707.69\\
G27 & 2000 & 3333   &  5.22  & {\textbf{3341}}   & 0.07  & {\textbf{3341}}  & 1.77 & {\textbf{3341}} & 585.01\\ 
G28 & 2000 & 3291   &  5.22  & {\textbf{3298}}   & 0.20  & 3287  & 1.77 & {\textbf{3298}} & 583.88\\ 
G29 & 2000 & 3395   &  5.22  & {\textbf{3405}}   & 0.07  & 3391  & 1.82 & {\textbf{3405}} & 572.47\\ 
G30 & 2000 & 3409   &  5.21  & {\textbf{3413}}   & 0.20  & 3412  & 1.79 & 3409 & 608.00\\
G31 & 2000 & 3305   &  5.23  & {\textbf{3310}}   & 0.79  & 3306  & 1.77 & 3309 & 577.13\\
G32 & 2000 & 1398   &  5.22  & {\textbf{1410}}   & 2.08  & 1400  & 1.57 & 1408 & 569.81 \\ 
G33 & 2000 & 1372   &  5.22  & {\textbf{1382}}   & 9.77  & 1374  & 1.56 & 1380 & 584.98\\
G34 & 2000 & 1376   &  5.22  & {\textbf{1384}}   & 0.78  & 1376  & 1.58 & 1382 & 595.74
\\\bottomrule
\end{tabular}
\end{adjustbox}
\label{table4}
\end{table}

\begin{table}[H]
\centering
\caption{A comparison between the best solutions obtained with the Simulated Annealing heuristic (SA), Toshiba Simulated  Bifurcation Machine solver (SBM), and the D-Wave Hybrid solver.}\label{tab:res4}
\begin{adjustbox}{width=1\textwidth}
\begin{tabular}{cccccccccc}
\toprule
Instance & $n$ & Hybrid & $t_{Hybrid} (s)$ & SBM & $t_{SBAM} (s)$  & SA1 & $t_{SA1} (s)$ & SA2 & $t_{SA2} (s)$ \\
\midrule
G35 & 2000 & 7659   &  5.22  & 7680   & 10.00 & 7472  & 2.41 & {\textbf{7686}} & 683.91\\
G36 & 2000 & 7655   &  5.21  & {\textbf{7675}}   & 10.01 & 7458  & 2.51 & {\textbf{7675}} & 668.88\\
G37 & 2000 & 7665   &  5.22  & {\textbf{7685}}   & 10.00 & 7464  & 2.41 & 7676 & 663.61\\
G38 & 2000 & 7661   &  5.22  & {\textbf{7686}}   & 10.00 & 7451  & 2.45 & 7680 & 689.77 \\ 
G39 & 2000 & 2390   &  5.22  & {\textbf{2407}}   & 10.00 & 2392  & 1.72 & 2400 & 580.16\\ 
G40 & 2000 & 2387   &  5.21  & {\textbf{2400}}   & 4.95  & 2386  & 1.70 & 2393 & 620.22 \\ 
G41 & 2000 & 2386   &  5.22  & 2404   & 10.00 & 2387  & 1.72 & {\textbf{2405}} & 592.46\\
G42 & 2000 & 2464   &  5.22  & {\textbf{2475}}   & 9.87  & 2465  & 1.73 & 2467 & 569.23\\
G43 & 1000 & {\textbf{6660}}   &  2.99  & {\textbf{6660}}   & 0.02  & 6656  & 2.13 & {\textbf{6660}} & 603.74\\
G44 & 1000 & {\textbf{6650}}   &  3.00  & {\textbf{6650}}   & 0.02  & 6637  & 2.14 & {\textbf{6650}} & 601.53\\ 
G45 & 1000 & 6653   &  2.99  & {\textbf{6654}}   & 0.02  & 6644  & 2.14 & {\textbf{6654}} & 601.96\\ 
G46 & 1000 & {\textbf{6649}}   &  2.99  & {\textbf{6649}}   & 0.02  & 6642  & 2.14 & {\textbf{6649}} & 604.65\\
G47 & 1000 & 6656   &  2.99  & {\textbf{6657}}   & 0.02  & 6647  & 2.13 & 6655 & 604.32\\ 
G48 & 3000 & {\textbf{6000}}   &  7.50  & {\textbf{6000}}   & 0.07  & {\textbf{6000}}  & 1.78 & {\textbf{6000}} & 558.03\\ 
G49 & 3000 & {\textbf{6000}}   &  7.50  & {\textbf{6000}}   & 0.07  & {\textbf{6000}}  & 1.75 & {\textbf{6000}} & 559.09\\ 
G50 & 3000 & {\textbf{5880}}   &  7.50  & {\textbf{5880}}   & 0.07  & 5824  & 1.78 & {\textbf{5880}} & 560.66\\ 
G51 & 1000 & 3843   &  2.99  & {\textbf{3848}}   & 0.94  & 3836  & 2.01 & {\textbf{3848}} & 597.82\\ 
G52 & 1000 & 3845   &  2.99  &  {\textbf{3851}}   & 1.77  & 3839  & 2.03 & 3849 & 596.27\\
G53 & 1000 & 3845   &  3.00  &  {\textbf{3849}}   & 10.00 & 3827  & 2.00 & 3848 & 598.09\\
G54 & 1000 & 3845   &  3.00  &  {\textbf{3851}}   & 10.00 & 3827  & 2.00 & 3845 & 617.68\\ 
G55 & 5000 & 10264  &  14.59 &  {\textbf{10289}}  & 10.00 & 9367  & 2.44 &  {\textbf{10289}} & 680.00\\ 
G56 & 5000 & 3982   &  14.58 &  {\textbf{4008}}   & 9.80  & 3912  & 1.98 & 4007 & 592.00\\ 
G57 & 5000 & 3456   &  14.59 & 3480   & 10.00 & 3436  & 1.84 &  {\textbf{3486}} & 582.95  \\
G58 & 5000 & 19221  &  14.59 & {\textbf{19257}}  & 10.01 & 17392 & 3.23 & 19118 & 858.10\\ 
G59 & 5000 & 6025   &  14.60 & {\textbf{6067}}   & 10.00 & 5892  & 2.28 & 6065 & 624.34\\
G60 & 7000 & 14142  &  24.76 & {\textbf{14168}}  & 10.00 & 12120 & 2.72 & 14114 & 733.45\\ 
G61 & 7000 & 5750   &  24.75 & 5777   & 10.01 & 5445  & 2.19 & {\textbf{5788}} & 616.86\\
G62 & 7000 & 4818   &  24.74 & 4844   & 10.01 & 4632  & 2.06 & {\textbf{4856}} & 603.25\\ 
G63 & 7000 & 26948  &  24.75 & {\textbf{26986}}  & 10.01 & 23751 & 3.47 & 26384 & 967.89\\ 
G64 & 7000 & 8665   &  24.74 & {\textbf{8728}}   & 8.82  & 8158  & 2.51 & {\textbf{8728}} & 661.92\\ 
G65 & 8000 & 5494   &  29.83 & 5532   & 10.01 & 5164  & 2.06 & {\textbf{5544}} & 613.44 \\
G66 & 9000 & 6282   &  34.92 & 6324   & 9.01  & 5770  & 2.13 & {\textbf{6350}} & 621.45\\ 
G67 & 10000 & 6868  &  39.98 & 6906   & 10.01 & 6136  & 2.18 & {\textbf{6924}} & 630.33\\
G70 & 10000 & 9516  &  40.00 & 9522   & 10.00 & 7909  & 2.48 & {\textbf{9546}} & 709.99\\
G72 & 10000 & 6914  &  39.98 & 6966   & 10.00 & 6148  & 2.20 & {\textbf{6984}} & 641.89
\\\bottomrule
\end{tabular}
\end{adjustbox}
\label{table4}
\end{table}

\begin{table}[H]
\centering
\caption{For each set of instances, we report how many times each algorithm discovered the best and worst solution. The "only" column counts the number of instances where the solver is the only one that achieved the best value, whereas the "worst" column gives the number of instances the solver computed the smallest value among all solvers, and this value was strictly smaller than the best value, although there may be another solver which also computed the same smallest value.}
\begin{adjustbox}{width=1\textwidth}
\begin{tabular}{cc|ccc|ccc|ccc|ccc}
\toprule
\multirow{2}{*}{Set} & \multirow{2}{*}{Set size} & \multicolumn{3}{c|}{Hybrid} & \multicolumn{3}{c|}{SBM} & \multicolumn{3}{c|}{SA1} & \multicolumn{3}{c}{SA2}\\
                     &                           &best & only & worst & best & only & worst & best& only & worst & best& only & worst\\
\midrule
{\tt be}  & 50 & 50 & 0 & 0 & -&- &- & 50 & 0 & 0 & 50&0 & 0  \\
{\tt bqp} & 20 & 20 & 0 & 0 & - & - & - & 20 & 0 & 0 & 20 & 0 & 0 \\
{\tt G}   & 69 & 17 & 0 & 11 & 58 & 26 & 0  & 8 & 0 & 53  & 41 & 11 & 3
\end{tabular}
\end{adjustbox}
\label{tab:res5}
\end{table}

\section{Discussion and conclusions}
In this paper, we presented a comprehensive analysis of exact, quantum, and classical approximate solvers for the NP-hard Max-Cut problem, focusing on the solution quality. The main objectives of this study were to evaluate the performance of the D-Wave fast annealing quantum processing unit (QPU) and Hybrid solvers in generating near-optimal solutions and to compare these results both with the best-known solutions from the literature and with solutions obtained using classical heuristic algorithms, in particular Simulated Annealing (algorithms SA1 and SA2) and the Toshiba Simulated Bifurcation Machine (SBM), for which we did not have the code available, but we referred to the results from the literature.

We utilized three benchmark datasets comprising 139 instances ranging from 100 to 10,000 nodes. For smaller instances (up to 251 nodes), we provided the global optima computed by the solvers BiqBin and MADAM. For larger instances, we used the best-known solutions from the literature. We have made these datasets publicly available for further numerical analysis.

On the first dataset, containing instances with a maximum size of 151 nodes, we were able to test all algorithms except SBM. Our results showed that the Hybrid solver and both simulated annealing variants consistently achieved the global optimum, whereas the QPU produced solutions that were far from optimal.

For the second dataset, the exact solver computed the global optimum for instances up to 251 nodes and matched the best-known values from the literature for instances of size 501. On this dataset, we tested the exact solver, the Hybrid solver, and both simulated annealing variants. The instances were too large for the QPU to handle. 
We observed similar results as on the first dataset, except the QPU could no longer be tested on these larger instances. The Hybrid solver and both simulated annealing variants achieved solutions with objective values that matched global optima where known and best-known values otherwise.

The third dataset, containing instances between 800 and 10,000 nodes, was especially interesting. For these instances, we referenced and reported the best-known solutions from the literature. Since these instances were too large for the exact solver and QPU, we applied only the Hybrid solver, both simulated annealing variants (SA1 and SA2), and reported SBM results from~\cite{matsuda2019benchmarking}. In addition to reporting Max-Cut solution values, we also included computational times (in seconds).
Our findings on this dataset showed that SBM and SA2 were highly competitive in solution quality, though SA2 required significantly more computation time. The Hybrid solver provided slightly lower-quality solutions compared to SA2 but achieved computation times closer to those of SBM. SA1, while comparable in computation time to the Hybrid solver and sometimes slightly slower than SBM, generally produced solutions of lower quality than the other three solvers.

Our results indicate that the D-Wave fast annealing quantum processing unit (QPU) serves as an intriguing solver, capable of producing multiple good solutions within a short time frame. However, these solutions are far from optimal, and the QPU's limitations in handling larger problem instances are significant; we could only solve instances up to a size of 151.

In contrast, the D-Wave Hybrid solver is highly efficient: it can handle large instances (up to 10,000) very quickly, producing solutions that are almost always optimal. For scientific research, however, it is frustrating that the Hybrid solver operates as a black-box cloud solution, leaving us without insight into its internal workings or the conditions under which it utilizes quantum resources.

The SBM also seems to be a highly competitive solver. Unfortunately, we had no access to this solver and were unable to conduct testing with it and had to rely entirely on publication~\cite{matsuda2019benchmarking}, which was not even published in a scientific journal, limiting our ability to fully evaluate this solver.

Notably, classical simulated annealing algorithms, which are transparent and can leverage specific features of the Max-Cut problem, compete strongly with the Hybrid solver in terms of both solution quality and time. The key advantage here is that we fully understand the steps within these classical algorithms.

Our study, therefore, advances the understanding of quantum-assisted optimization for combinatorial challenges and offers valuable insights for selecting appropriate solvers based on application-specific constraints and desired performance.

\begin{acknowledgements}
The research of J.P. was co-funded by the Republic of Slovenia, the Ministry of HE, Science and Innovation, and the Slovenian Research and Innovation Agency of the European Union—NextGenerationEU through the DIGITOP project,  by the Slovenian Research
and Innovation Agency (ARIS) through the annual work program of Rudolfovo and the research project L1-60136.

J.V. acknowledges the financial support from ARIS, P08333 young researcher grant, and P1-0040 Nonequilibrium Quantum System Dynamics.

\end{acknowledgements}

% Non-BibTeX users please use
%\bibliography{References.bib}

\printbibliography

@article{dunning2018works,
  title={What works best when? A systematic evaluation of heuristics for Max-Cut and QUBO},
  author={Dunning, Iain and Gupta, Swati and Silberholz, John},
  journal={INFORMS Journal on Computing},
  volume={30},
  number={3},
  pages={608--624},
  year={2018},
  publisher={INFORMS}
}

@misc{matsuda2019benchmarking,
  title={Benchmarking the MAX-CUT problem on the Simulated Bifurcation Machine},
  author={Matsuda, Y},
  year={2019},
  publisher={Medium},
url={https://medium.com/toshiba-sbm/benchmarking-the-max-cut-problem-on-the-simulated-bifurcation-machine-e26e1127c0b0}
}

@misc{biqmaclib,
author={Wiegele, Angelika},
title={Biq{M}ac {L}ibrary},
year=2007,
howpublished={\url{http://biqmac.aau.at/biqmaclib.html}}
}

@article{wiegele2007biq,
  title={Biq Mac Library -- A collection of Max-Cut and quadratic 0-1 programming instances of medium size},
  author={Wiegele, Angelika},
  journal={Preprint},
  year={2007}
}

@article{kirkpatrick1983optimization,
  title={Optimization by simulated annealing},
  author={Kirkpatrick, Scott and Gelatt Jr, C Daniel and Vecchi, Mario P},
  journal={science},
  volume={220},
  number={4598},
  pages={671--680},
  year={1983},
  publisher={American association for the advancement of science}
}

@article{vcerny1985thermodynamical,
  title={Thermodynamical approach to the traveling salesman problem: An efficient simulation algorithm},
  author={{\v{C}}ern{\`y}, Vladim{\'\i}r},
  journal={Journal of optimization theory and applications},
  volume={45},
  pages={41--51},
  year={1985},
  publisher={Springer}
}

@article{goemans1995improved,
  title={Improved approximation algorithms for maximum cut and satisfiability problems using semidefinite programming},
  author={Goemans, Michel X and Williamson, David P},
  journal={Journal of the ACM (JACM)},
  volume={42},
  number={6},
  pages={1115--1145},
  year={1995},
  publisher={ACM New York, NY, USA}
}

@article{khot2007optimal,
  title={Optimal inapproximability results for MAX-CUT and other 2-variable CSPs?},
  author={Khot, Subhash and Kindler, Guy and Mossel, Elchanan and O’Donnell, Ryan},
  journal={SIAM Journal on Computing},
  volume={37},
  number={1},
  pages={319--357},
  year={2007},
  publisher={SIAM}
}

@article{myklebust2015solving,
  title={Solving maximum cut problems by simulated annealing},
  author={Myklebust, Tor GJ},
  journal={arXiv preprint arXiv:1505.03068},
  year={2015}
}

@article{helmberg2000spectral,
  title={A spectral bundle method for semidefinite programming},
  author={Helmberg, Christoph and Rendl, Franz},
  journal={SIAM Journal on Optimization},
  volume={10},
  number={3},
  pages={673--696},
  year={2000},
  publisher={SIAM}
}

@article{willsch2022benchmarking,
  title={Benchmarking Advantage and D-Wave 2000Q quantum annealers with exact cover problems},
  author={Willsch, Dennis and Willsch, Madita and Gonzalez Calaza, Carlos D and Jin, Fengping and De Raedt, Hans and Svensson, Marika and Michielsen, Kristel},
  journal={Quantum Information Processing},
  volume={21},
  number={4},
  pages={141},
  year={2022},
  publisher={Springer}
}

@article{sharma2023augmenting,
  title={Augmenting an electronic Ising machine to effectively solve boolean satisfiability},
  author={Sharma, Anshujit and Burns, Matthew and Hahn, Andrew and Huang, Michael},
  journal={Scientific Reports},
  volume={13},
  number={1},
  pages={22858},
  year={2023},
  publisher={Nature Publishing Group UK London}
}

@inproceedings{montanez2023improving,
  title={Improving Performance in Combinatorial Optimization Problems with Inequality Constraints: An Evaluation of the Unbalanced Penalization Method on D-Wave Advantage},
  author={Montanez-Barrera, JA and van den Heuvel, Pim and Willsch, Dennis and Michielsen, Kristel},
  booktitle={2023 IEEE International Conference on Quantum Computing and Engineering (QCE)},
  volume={1},
  pages={535--542},
  year={2023},
  organization={IEEE}
}

@article{kowalsky20223,
  title={3-regular three-XORSAT planted solutions benchmark of classical and quantum heuristic optimizers},
  author={Kowalsky, Matthew and Albash, Tameem and Hen, Itay and Lidar, Daniel A},
  journal={Quantum Science and Technology},
  volume={7},
  number={2},
  pages={025008},
  year={2022},
  publisher={IOP Publishing}
}

@article{tasseff2022emerging,
  title={On the emerging potential of quantum annealing hardware for combinatorial optimization},
  author={Tasseff, Byron and Albash, Tameem and Morrell, Zachary and Vuffray, Marc and Lokhov, Andrey Y and Misra, Sidhant and Coffrin, Carleton},
  journal={arXiv preprint arXiv:2210.04291},
  year={2022}
}

@inproceedings{mcleod2022benchmarking,
  title={Benchmarking d-wave quantum annealers: Spectral gap scaling of maximum cardinality matching problems},
  author={McLeod, Cameron Robert and Sasdelli, Michele},
  booktitle={International Conference on Computational Science},
  pages={150--163},
  year={2022},
  organization={Springer}
}

@article{konz2021embedding,
  title={Embedding overhead scaling of optimization problems in quantum annealing},
  author={K{\"o}nz, Mario S and Lechner, Wolfgang and Katzgraber, Helmut G and Troyer, Matthias},
  journal={PRX Quantum},
  volume={2},
  number={4},
  pages={040322},
  year={2021},
  publisher={APS}
}

@article{koshikawa2021benchmark,
  title={Benchmark test of black-box optimization using d-wave quantum annealer},
  author={Koshikawa, Ami S and Ohzeki, Masayuki and Kadowaki, Tadashi and Tanaka, Kazuyuki},
  journal={Journal of the Physical Society of Japan},
  volume={90},
  number={6},
  pages={064001},
  year={2021},
  publisher={The Physical Society of Japan}
}

@article{gonzalez2021garden,
  title={Garden optimization problems for benchmarking quantum annealers},
  author={Gonzalez Calaza, Carlos D and Willsch, Dennis and Michielsen, Kristel},
  journal={Quantum Information Processing},
  volume={20},
  number={9},
  pages={305},
  year={2021},
  publisher={Springer}
}

@online{dwave2021performanceupdate,
  author = {D-Wave Systems Inc},
  title = {The Advantage System: Performance Update},
  year = 2021,
  url = {https://www.dwavesys.com/media/kjtlcemb/14-1054a-a_advantage_system_performance_update.pdf},
  urldate = {2024-04-22}
}

@article{zlokapa2020quantum,
  title={Quantum adiabatic machine learning by zooming into a region of the energy surface},
  author={Zlokapa, Alexander and Mott, Alex and Job, Joshua and Vlimant, Jean-Roch and Lidar, Daniel and Spiropulu, Maria},
  journal={Physical Review A},
  volume={102},
  number={6},
  pages={062405},
  year={2020},
  publisher={APS}
}

@article{pang2021potential,
  title={The potential of quantum annealing for rapid solution structure identification},
  author={Pang, Yuchen and Coffrin, Carleton and Lokhov, Andrey Y and Vuffray, Marc},
  journal={Constraints},
  volume={26},
  number={1},
  pages={1--25},
  year={2021},
  publisher={Springer}
}

@article{perdomo2019readiness,
  title={Readiness of quantum optimization machines for industrial applications},
  author={Perdomo-Ortiz, Alejandro and Feldman, Alexander and Ozaeta, Asier and Isakov, Sergei V and Zhu, Zheng and O’Gorman, Bryan and Katzgraber, Helmut G and Diedrich, Alexander and Neven, Hartmut and de Kleer, Johan and others},
  journal={Physical Review Applied},
  volume={12},
  number={1},
  pages={014004},
  year={2019},
  publisher={APS}
}

@article{king2019quantum,
  title={Quantum annealing amid local ruggedness and global frustration},
  author={King, James and Yarkoni, Sheir and Raymond, Jack and Ozfidan, Isil and King, Andrew D and Nevisi, Mayssam Mohammadi and Hilton, Jeremy P and McGeoch, Catherine C},
  journal={Journal of the Physical Society of Japan},
  volume={88},
  number={6},
  pages={061007},
  year={2019},
  publisher={The Physical Society of Japan}
}

@article{hamerly2019experimental,
  title={Experimental investigation of performance differences between coherent Ising machines and a quantum annealer},
  author={Hamerly, Ryan and Inagaki, Takahiro and McMahon, Peter L and Venturelli, Davide and Marandi, Alireza and Onodera, Tatsuhiro and Ng, Edwin and Langrock, Carsten and Inaba, Kensuke and Honjo, Toshimori and others},
  journal={Science advances},
  volume={5},
  number={5},
  pages={eaau0823},
  year={2019},
  publisher={American Association for the Advancement of Science}
}

@article{zhu2019fair,
  title={Fair sampling of ground-state configurations of binary optimization problems},
  author={Zhu, Zheng and Ochoa, Andrew J and Katzgraber, Helmut G},
  journal={Physical Review E},
  volume={99},
  number={6},
  pages={063314},
  year={2019},
  publisher={APS}
}

@article{katzgraber2018viewing,
  title={Viewing vanilla quantum annealing through spin glasses},
  author={Katzgraber, Helmut G},
  journal={Quantum Science and Technology},
  volume={3},
  number={3},
  pages={030505},
  year={2018},
  publisher={IOP Publishing}
}

@article{albash2017temperature,
  title={Temperature scaling law for quantum annealing optimizers},
  author={Albash, Tameem and Martin-Mayor, Victor and Hen, Itay},
  journal={Physical review letters},
  volume={119},
  number={11},
  pages={110502},
  year={2017},
  publisher={APS}
}

@article{mandra2017exponentially,
  title={Exponentially biased ground-state sampling of quantum annealing machines with transverse-field driving hamiltonians},
  author={Mandra, Salvatore and Zhu, Zheng and Katzgraber, Helmut G},
  journal={Physical review letters},
  volume={118},
  number={7},
  pages={070502},
  year={2017},
  publisher={APS}
}

@article{denchev2016computational,
  title={What is the computational value of finite-range tunneling?},
  author={Denchev, Vasil S and Boixo, Sergio and Isakov, Sergei V and Ding, Nan and Babbush, Ryan and Smelyanskiy, Vadim and Martinis, John and Neven, Hartmut},
  journal={Physical Review X},
  volume={6},
  number={3},
  pages={031015},
  year={2016},
  publisher={APS}
}

@article{mandra2016strengths,
  title={Strengths and weaknesses of weak-strong cluster problems: A detailed overview of state-of-the-art classical heuristics versus quantum approaches},
  author={Mandra, Salvatore and Zhu, Zheng and Wang, Wenlong and Perdomo-Ortiz, Alejandro and Katzgraber, Helmut G},
  journal={Physical Review A},
  volume={94},
  number={2},
  pages={022337},
  year={2016},
  publisher={APS}
}

@article{zhu2016best,
  title={Best-case performance of quantum annealers on native spin-glass benchmarks: How chaos can affect success probabilities},
  author={Zhu, Zheng and Ochoa, Andrew J and Schnabel, Stefan and Hamze, Firas and Katzgraber, Helmut G},
  journal={Physical Review A},
  volume={93},
  number={1},
  pages={012317},
  year={2016},
  publisher={APS}
}

@article{boixo2016computational,
  title={Computational multiqubit tunnelling in programmable quantum annealers},
  author={Boixo, Sergio and Smelyanskiy, Vadim N and Shabani, Alireza and Isakov, Sergei V and Dykman, Mark and Denchev, Vasil S and Amin, Mohammad H and Smirnov, Anatoly Yu and Mohseni, Masoud and Neven, Hartmut},
  journal={Nature communications},
  volume={7},
  number={1},
  pages={10327},
  year={2016},
  publisher={Nature Publishing Group UK London}
}

@article{cao2016solving,
  title={Solving set cover with pairs problem using quantum annealing},
  author={Cao, Yudong and Jiang, Shuxian and Perouli, Debbie and Kais, Sabre},
  journal={Scientific reports},
  volume={6},
  number={1},
  pages={33957},
  year={2016},
  publisher={Nature Publishing Group UK London}
}

@article{radovic2018machine,
  title={Machine learning at the energy and intensity frontiers of particle physics},
  author={Radovic, Alexander and Williams, Mike and Rousseau, David and Kagan, Michael and Bonacorsi, Daniele and Himmel, Alexander and Aurisano, Adam and Terao, Kazuhiro and Wongjirad, Taritree},
  journal={Nature},
  volume={560},
  number={7716},
  pages={41--48},
  year={2018},
  publisher={Nature Publishing Group UK London}
}

@article{king2024computational,
  title={Computational supremacy in quantum simulation},
  author={King, Andrew D and Nocera, Alberto and Rams, Marek M and Dziarmaga, Jacek and Wiersema, Roeland and Bernoudy, William and Raymond, Jack and Kaushal, Nitin and Heinsdorf, Niclas and Harris, Richard and others},
  journal={arXiv preprint arXiv:2403.00910},
  year={2024}
}

@article{king2023quantum,
  title={Quantum critical dynamics in a 5,000-qubit programmable spin glass},
  author={King, Andrew D and Raymond, Jack and Lanting, Trevor and Harris, Richard and Zucca, Alex and Altomare, Fabio and Berkley, Andrew J and Boothby, Kelly and Ejtemaee, Sara and Enderud, Colin and others},
  journal={Nature},
  volume={617},
  number={7959},
  pages={61--66},
  year={2023},
  publisher={Nature Publishing Group UK London}
}

@article{king2022coherent,
  title={Coherent quantum annealing in a programmable 2,000 qubit Ising chain},
  author={King, Andrew D and Suzuki, Sei and Raymond, Jack and Zucca, Alex and Lanting, Trevor and Altomare, Fabio and Berkley, Andrew J and Ejtemaee, Sara and Hoskinson, Emile and Huang, Shuiyuan and others},
  journal={Nature Physics},
  volume={18},
  number={11},
  pages={1324--1328},
  year={2022},
  publisher={Nature Publishing Group UK London}
}

@article{gusmeroli2022biqbin,
  title={BiqBin: a parallel branch-and-bound solver for binary quadratic problems with linear constraints},
  author={Gusmeroli, Nicolo and Hrga, Timotej and Lu{\v{z}}ar, Borut and Povh, Janez and Siebenhofer, Melanie and Wiegele, Angelika},
  journal={ACM Transactions on Mathematical Software (TOMS)},
  volume={48},
  number={2},
  pages={1--31},
  year={2022},
  publisher={ACM New York, NY}
}

@article{hrga2021madam,
  title={MADAM: a parallel exact solver for max-cut based on semidefinite programming and ADMM},
  author={Hrga, Timotej and Povh, Janez},
  journal={Computational Optimization and Applications},
  volume={80},
  number={2},
  pages={347--375},
  year={2021},
  publisher={Springer}
}

@article{arute2019quantum,
  title={Quantum supremacy using a programmable superconducting processor},
  author={Arute, Frank and Arya, Kunal and Babbush, Ryan and Bacon, Dave and Bardin, Joseph C and Barends, Rami and Biswas, Rupak and Boixo, Sergio and Brandao, Fernando GSL and Buell, David A and others},
  journal={Nature},
  volume={574},
  number={7779},
  pages={505--510},
  year={2019},
  publisher={Nature Publishing Group}
}

@article{hauke2020perspectives,
  title={Perspectives of quantum annealing: Methods and implementations},
  author={Hauke, Philipp and Katzgraber, Helmut G and Lechner, Wolfgang and Nishimori, Hidetoshi and Oliver, William D},
  journal={Reports on Progress in Physics},
  volume={83},
  number={5},
  pages={054401},
  year={2020},
  publisher={IOP Publishing}
}

@article{barahona1989experiments,
  title={Experiments in quadratic 0--1 programming},
  author={Barahona, Francisco and J{\"u}nger, Michael and Reinelt, Gerhard},
  journal={Mathematical programming},
  volume={44},
  number={1},
  pages={127--137},
  year={1989},
  publisher={Springer}
}

@book{karp2010reducibility,
  title={Reducibility among combinatorial problems},
  author={Karp, Richard M},
  year={2010},
  publisher={Springer}
}

@article{barahona1988application,
  title={An application of combinatorial optimization to statistical physics and circuit layout design},
  author={Barahona, Francisco and Gr{\"o}tschel, Martin and J{\"u}nger, Michael and Reinelt, Gerhard},
  journal={Operations Research},
  volume={36},
  number={3},
  pages={493--513},
  year={1988},
  publisher={INFORMS}
}

@phdthesis{liers2004contributions,
  title={Contributions to determining exact ground-states of Ising spin-glasses and to their physics},
  author={Liers, Frauke},
  year={2004},
  school={Universit{\"a}t zu K{\"o}ln}
}

@inproceedings{grotschel1987calculating,
  title={Calculating exact ground states of spin glasses: A polyhedral approach},
  author={Gr{\"o}tschel, Martin and J{\"u}nger, Michael and Reinelt, Gerhard},
  booktitle={Heidelberg Colloquium on Glassy Dynamics: Proceedings of a Colloquium on Spin Glasses, Optimization and Neural Networks Held at the University of Heidelberg June 9--13, 1986},
  pages={325--353},
  year={1987},
  organization={Springer}
}

@article{junger2021quantum,
  title={Quantum annealing versus digital computing: An experimental comparison},
  author={J{\"u}nger, Michael and Lobe, Elisabeth and Mutzel, Petra and Reinelt, Gerhard and Rendl, Franz and Rinaldi, Giovanni and Stollenwerk, Tobias},
  journal={Journal of Experimental Algorithmics (JEA)},
  volume={26},
  pages={1--30},
  year={2021},
  publisher={ACM New York, NY, USA}
}

@article{albash2018demonstration,
  title={Demonstration of a scaling advantage for a quantum annealer over simulated annealing},
  author={Albash, Tameem and Lidar, Daniel A},
  journal={Physical Review X},
  volume={8},
  number={3},
  pages={031016},
  year={2018},
  publisher={APS}
}

@article{king2021scaling,
  title={Scaling advantage over path-integral Monte Carlo in quantum simulation of geometrically frustrated magnets},
  author={King, Andrew D and Raymond, Jack and Lanting, Trevor and Isakov, Sergei V and Mohseni, Masoud and Poulin-Lamarre, Gabriel and Ejtemaee, Sara and Bernoudy, William and Ozfidan, Isil and Smirnov, Anatoly Yu and others},
  journal={Nature communications},
  volume={12},
  number={1},
  pages={1113},
  year={2021},
  publisher={Nature Publishing Group UK London}
}

@article{pelofske2024increasing,
  title={Increasing the Hardness of Posiform Planting Using Random QUBOs for Programmable Quantum Annealer Benchmarking},
  author={Pelofske, Elijah and Hahn, Georg and Djidjev, Hristo},
  journal={arXiv preprint arXiv:2411.03626},
  year={2024}
}

@software{erzen_2024_14290290,
  author       = {Eržen, Vid and
                  Hrga, Timotej and
                  Povh, Janez and
                  Vodeb, Jaka},
  title        = {Max-Cut benchmark dataset},
  month        = dec,
  year         = 2024,
  publisher    = {Zenodo},
  version      = {v0.0.4},
  doi          = {10.5281/zenodo.14290290},
  url          = {https://doi.org/10.5281/zenodo.14290290}
 }

\end{document}